\documentclass[10pt]{paper}

\usepackage{amssymb,amsmath,enumerate,theorem,epsfig,rotating,graphics,changebar,eepic}
\usepackage{amsfonts}
\usepackage{a4,latexsym,parskip}
\usepackage{hyperref}
\hypersetup{pdfauthor=Matthias Sch"utt} \hypersetup{pdftitle=Diss}

\newcommand{\C} [1][]{\mathbb{C}^{#1}}
\newcommand{\Q} [1] []{\mathbb{Q}_{#1}}
\newcommand{\N} [1][] {\mathbb{N}_{#1}}

\newcommand{\F}{\mathbb{F}}
\newcommand{\Z}{\mathbb{Z}}
\newcommand{\p}{\mathfrak{p}}
\newcommand{\m}{\mathfrak{m}}
\newcommand{\OO}{\mathcal{O}}
\newcommand{\NN}{\mathcal{N}}
\newcommand{\CS}{\mathcal{S}}

\newcommand{\e}{\varepsilon}
\newcommand{\ez}{\eta_{\Z}}
\newcommand{\ia}{\mathfrak{a}}
\newcommand{\jj}{\mathfrak{j}}
\newcommand{\teta}{\tilde\eta}

\newcommand{\qed}{\hfill \ensuremath{\Box}}

\theoremstyle{break} \newtheorem{Theorem}{Theorem}[section]
\newtheorem{Proposition}[Theorem]{Proposition}
\newtheorem{Lemma}[Theorem]{Lemma}

\newtheorem{Definition}[Theorem]{Definition}
\newtheorem{Corollary}[Theorem]{Corollary}

\newtheorem{Remark}[Theorem]{Remark}
\newtheorem{Question}[Theorem]{Question}

\begin{document}
\setlength{\unitlength}{1cm}

\title{CM newforms with rational coefficients}
\author{Matthias Sch\"utt}
\maketitle

\abstract{We classify newforms with rational Fourier coefficients and complex
multiplication for fixed weight up to twisting. 
Under the extended Riemann hypothesis for odd real Dirichlet characters,
these newforms are finite in number.
We produce tables for weights 3 and 4, where finiteness holds
unconditionally.}

\begin{small}

\end{small}

\vspace{0.2cm}

\section*{Introduction}
\label{intro}

This paper grew out of the wish to understand and eventually
determine the modular forms which could possibly be associated to
singular K3 surfaces over $\Q$. The term "singular" refers to the
maximal Picard number $\varrho=20$. Hence the transcendental
lattice $T_X$ of a singular K3 surface $X$ has rank two, giving
rise to a two-dimensional Galois representation. By a theorem of
Livn\'e \cite{L}, any singular K3 surface over $\Q$ is modular,
with the associated modular form of weight 3. Then we can ask
which newforms occur. Our original motivation was to determine all
newforms of weight 3 with rational coefficients, that is, all
possibly occurring newforms. Up to twisting, we will produce a
finite list in Table \ref{T:wt3}.

All the ideas involved in our analysis can be applied to newforms
of general weight which have rational coefficients and CM. Hence
we will treat this question in the most possible generality. A
guiding and illustrating example is provided by the newforms of
weight 2 with rational Fourier coefficients and CM. Since these
correspond to elliptic curves over $\Q$ with CM, the finiteness
result is guaranteed by the classical theory of elliptic curves
with CM, and in particular, by class field theory. Our main result is as follows:

\textbf{Theorem \ref{Thm:mod}}\\
{\sl
Assume the extended Riemann hypothesis (ERH) for odd real Dirichlet characters. Then there are only finitely many CM newforms with
rational coefficients for fixed weight up to twisting.
For weights $2,3,4,5,6, 2^r+1, 3\cdot 2^r+1 (r\in\N)$ this holds unconditionally.}

The paper is organised as follows: First we recall some of the
basic theory of modular forms and Hecke characters (Sect.~\ref{s:mf}). Then we
state our results concerning the finiteness and shape of
CM newforms with rational coefficients for fixed weight.
Section \ref{s:results} also explains the dependence of Theorem \ref{Thm:mod} 
on the ERH.
The subsequent sections give a
proof of the results. We conclude the
paper by an explicit description of CM newforms of weight 3 and 4
with rational coefficients (Sect.~\ref{s:wt3}, \ref{s:wt4}). In
particular, this will answer the motivating question. Finally, we comment on 
recent progress on the problem of
geometric realisations of CM newforms with rational coefficients.
Here we also return to singular K3 surfaces as Thm.~\ref{Thm:mod} enables N.~Elkies and the author to solve the geometric realisation problem in weight 3 (cf.~Thm.~\ref{Thm:ES}).

\section{Modular forms and Hecke characters}
\label{s:mf}
\label{s:Gr}

In this section we recall some facts about modular forms and sketch the connection between eigenforms with complex
multiplication (CM) and Hecke characters. The fact
that a Hecke character gives rise to a modular form goes
back to Hecke. On the other hand, Ribet showed that any newform
with CM comes from a Hecke character (Thm.~\ref{Thm:Ribet}).
For details, the reader is referred to the paper of Ribet
\cite{R}.

We are concerned with newforms in the space of (elliptic) cusp
forms $\CS_k(\Gamma_0(N),\e)$. Here, we fix an integral weight $k>1$, the level $N\in\N$ and the \emph{nebentypus}, a character $\e$ modulo $N$ with $\e(-1)=(-1)^k$. A cusp form has a Fourier expansion at $i\infty$ 
\[
f(\tau)=\sum_{n=1}^\infty a_nq^n,\;\;\;\;\; \text{ where } \;\,q=e^{2\pi i\tau}.
\]
A normalised newform $f$ can be
characterised by the minimality of its level and the property that
its Mellin transform $L(f,s)$ possesses an Euler product expansion
\[
L(f,s)=\sum_{n\in\N} a_n n^{-s}=\prod_p (1-a_pp^{-s}+\e(p)p^{k-1-2s})^{-1}.
\]
Consider the finite
extension of $\Q$ which is generated by the set of Fourier
coefficients $\{a_p\}$ of a newform $f$. This field is well-known to be
real if and only if the nebentypus $\e$ is either trivial
or quadratic with $\e(p)a_p=a_p$; otherwise it is a CM-field.

A newform $f=\sum a_n q^n\in\CS_k(\Gamma_1(N)$ is said to have
\emph{complex multiplication (CM)} by a Dirichlet character $\phi$
if $f=f\otimes\phi$. Here, we write the twist
\[
f\otimes\phi=\sum_{n\in\N} \phi(n) a_n q^n.
\]
Since a CM-character is neccessarily quadratic, we will also refer to CM by the corresponding quadratic field $K$. The ring of integers in $K$ will be denoted by $\OO_K$.

\begin{Theorem}[Ribet {\cite[Prop.~(4.4), Thm.~(4.5)]{R}}]\label{Thm:Ribet}
A newform has CM by $K$ (a quadratic field) if and only if it comes from a Hecke character of $K$.
In particular, the field $K$ is imaginary and unique.
\end{Theorem}

\begin{Corollary}\label{Cor:coeff}
Let $f$ be a newform with real coefficients. Then the following relations hold:
\begin{enumerate}[(i)]
\item If the weight $k$ is even, then the nebentypus $\e$ is
trivial. \item If the weight $k$ is odd, then $f$ has CM by its
nebentypus $\e$. In particular, $\e$ is quadratic.
\end{enumerate}
\end{Corollary}

\begin{Definition}[Hecke]
Let $\m$ be an ideal of $K$ and $l\in\N$. A
\emph{Hecke character} $\psi$ of $K$ modulo $\mathfrak{m}$
with $\infty$-type $l$ is a homomorphism on the group of
fractional ideals of $K$ which are prime to $\mathfrak{m}$,
\[
\psi: \left\{\begin{matrix}
\text{fractional ideals of $K$} \\ \text{ which are prime to } \mathfrak{m} \end{matrix} \right\} \rightarrow \C[*],
\]
such that
\[
\psi(\alpha\,\OO_K)=\alpha^l\;\;\;\;\;\; \forall\, \alpha\in K^* \text{ with }\alpha\equiv 1\mod\mathfrak{m}.
\]
The ideal $\m$ is called the \emph{conductor} of $\psi$ if $\m$ is
minimal in the following sense: If $\psi$ is defined modulo $\m'$,
then $\m|\m'$.
\end{Definition}
Strictly speaking, we should refer to multiplicative congruence in
this context. To simplify notation, we shall omit this in the
following and also tacitly extend $\psi$ by 0 for all fractional
ideals of $K$ which are not prime to $\m$.

The connection with modular forms occurs through the $L$-series of $\psi$. Denote the norm homomorphism of the Galois extension $K/\Q$ by $\mathcal{N}$:
\[
L(\psi,s)=\sum_{\mathfrak{a} \text{ integral}}\psi(\mathfrak{a})\, 
\mathcal{N}(\mathfrak{a})^{-s}=\prod_{\mathfrak{p} \text{ prime}}(1-\psi(\mathfrak{p}) 
\,\mathcal{N}(\mathfrak{p})^{-s})^{-1}.
\]
Then the inverse Mellin transform $f_\psi$ of $L(\psi,s)$ is known to be a Hecke eigenform:
\[
f_\psi=\sum_{n\in\N}a_nq^n=\sum_{\mathfrak{a} \text{ integral}}\psi(\mathfrak{a}) q^{\mathcal{N}(\mathfrak{a})} 
\;\;\;\;\;\;\;\;\;(q=e^{2\pi i\tau}).
\]
We introduce the following notation: Let $-\Delta_K$ be the discriminant of $K$ and $\chi_K$ the associated quadratic character of conductor $\Delta_K$. Denote $M=\mathcal{N}(\mathfrak{m})$ and consider the 
Dirichlet character $\eta_{\Z}$ mod $M$ which is given by
\[
\eta_{\Z}: a \mapsto \frac{\psi(a\,\OO_K)}{a^l} \;\;\;\;  (a\in\Z, (a,M)=1).
\]
\begin{Theorem}[Hecke, Shimura]\label{Thm:Shimura}
$f_\psi$ is a Hecke eigenform of weight $l+1$, level $\Delta_K M$ and
nebentypus character $\chi_K \eta_{\Z}:$
\[
f_\psi\in \CS_{l+1}(\Gamma_0(\Delta_K M),\chi_K \eta_{\Z}).
\]
$f_\psi$ is a newform (of level $\Delta_K M$) if and only if $\m$ is
the conductor of $\psi$.
\end{Theorem}

By construction, the eigenform $f_\psi$ has CM by $K$ (or
equivalently $\chi_K$). From Corollary \ref{Cor:coeff}, we obtain
the following restrictions on the nebentypus $\e=\chi_K\eta_{\Z}$:
\begin{Corollary}\label{Cor:ez}
Assume that the newform $f_\psi$ has real coefficients.
Then the following conditions hold:
\begin{enumerate}[(i)]
\item If $l$ is odd, then $\e=1$ and hence $\eta_{\Z}=\chi_K$. In particular, $(\sqrt{-\Delta_K})\mid\m$.
\item If $l$ is even, then $\e=\chi_K$ and hence $\eta_{\Z}=1$.
\end{enumerate}
\end{Corollary}

We say that a Hecke character $\psi$ has real (or rational) coefficients, if the corresponding eigenform $f_\psi$ has. 

\begin{Remark}\label{Rem:ez}
The first statement of Corollary \ref{Cor:ez} can be strengthened in two special cases as follows: If $\psi$ is a Hecke character of $\Q(\sqrt{-1})$ or $\Q(\sqrt{-2})$ with real coefficients and odd $\infty$-type, then
\[
(2+2\sqrt{-1})\mid\m \;\;\text{ resp.~}\;\; (4\sqrt{-2})\mid\m.
\]
\end{Remark}

If we twist $f_\psi$ by a Dirichlet character $\phi$, this corresponds to twisting $\psi$ by $\phi\circ\NN$:
\[
f_{\psi}\otimes\phi=f_{\psi\otimes(\phi\circ\NN)}.
\]
In general, the coefficients stay real only if $\phi$ is quadratic. If $K=\Q(\sqrt{-1})$ or $\Q(\sqrt{-3})$, we can twist $\psi$ also by quartic resp.~sextic characters.

\pagebreak

\section{Formulation of the results}
\label{s:results}

In this section we formulate our main result for CM newforms
(Thm.~\ref{Thm:mod}). For the proof we translate the statement to
Hecke characters (Thm.~\ref{Thm:Gr}). Then Thm.~\ref{Thm:Gr} follows from an unconditional classification of
Hecke characters (Thm.~\ref{1:1}). 
The conditional part of Thm.~\ref{Thm:mod} and \ref{Thm:Gr} 
is based on work of Weinberger \cite{Wb}, as we explain
in the sequel of Thm.~\ref{1:1}.
The unconditional part also relies on work by Heath-Brown. The proof of Thm.~\ref{1:1}
will be the subject of the next three sections.

\begin{Theorem}\label{Thm:mod}
Assume ERH for odd real Dirichlet characters. Then there are only finitely many CM newforms with
rational coefficients for fixed weight up to twisting.\\
For weights $2,3,4,5,6, 2^r+1, 3\cdot 2^r+1 (r\in\N)$ this holds unconditionally.
\end{Theorem}

For weight 2 this result is known from the theory of elliptic
curves with CM. For the general case, we will prove a 
corresponding statement for Hecke characters in Theorem
\ref{Thm:Gr}. Then we can use Ribet's Theorem \ref{Thm:Ribet} to
deduce Theorem \ref{Thm:mod}. For this purpose, we need the
following easy lemma:

\begin{Lemma}\label{Lem:rat}
Let $\psi$ be a Hecke character of an imaginary quadratic
field $K$. Then the
following statements are equivalent:
\begin{enumerate}[(i )]
\item $f_\psi$ has rational coefficients,
\item $\psi(\ia)+\psi(\bar\ia)\in\Z\;\;$ for all ideals $\ia\subset \OO_K$. 
\item im $\psi\subseteq\OO_K$.
\end{enumerate}
\end{Lemma}

As a consequence, Theorem \ref{Thm:mod} is equivalent to

\begin{Theorem}\label{Thm:Gr}
Assume ERH for odd real Dirichlet characters. For fixed $\infty$-type there are only finitely many Hecke characters of quadratic imaginary fields with rational coefficients up to  twisting. 
For $\infty$-types $1, 2, 3, 4, 5, 2^r, 3\cdot 2^r (r\in\N)$ this holds unconditionally.
\end{Theorem}

To derive Thm.~\ref{Thm:Gr}, we need the following unconditional classification. It involves the \emph{exponent} $e_K$ of an imaginary quadratic field $K$, that
is, the exponent of the class group $Cl(\OO_K)$.

\begin{Theorem}\label{1:1}
For fixed $\infty$-type $l$, there is a bijective correspondence
\[
\left\{\begin{matrix}
\text{Hecke characters}\\ \text{with $\infty$-type $l$ and rational coefficients}\\ \text{up to  twisting}
       \end{matrix}\right\}
\stackrel{1:1}{\longleftrightarrow}\left\{\begin{matrix}
\text{Imaginary quadratic fields $K$}\\ \text{with exponent}\;\, e_K\mid l
                                          \end{matrix}\right\}.
\]
\end{Theorem}

Theorem \ref{1:1} implies Theorem
\ref{Thm:Gr} due to a result of Weinberger: Subject to a Siegel-Landau condition
on the zeroes of the Dirichlet $L$-functions $L(s,\chi_K)$ for
imaginary quadratic fields $K$, he shows that
\[
e_K\to\infty \;\;\text{ as }\; \Delta_K\to\infty
\;\;\;\;\text{\cite[Thm.~3]{Wb}}.
\]
In particular this holds if the extended Riemann hypothesis for
$L(s,\chi_K)$ is true (cf.~\cite[Thm.~4]{Wb}). Hence, Theorem
\ref{1:1} implies the conditional part of Theorem \ref{Thm:Gr} as
claimed. A formulation of Theorem \ref{Thm:Gr}
with the weaker condition is omitted here for brevity.

Furthermore, Weinberger proves the finiteness of imaginary
quadratic fields $K$ with exponent $e_K=2$ or $3$ unconditionally
\cite[Thms. 1, 2]{Wb} (the latter is also due to Boyd-Kisilevski
\cite{BK}). Finiteness for the greater exponents in the extra statement of Theorem \ref{Thm:Gr} was recently proven by Heath-Brown \cite{HB}.

We will pay special attention to unconditional cases in
sections \ref{s:wt3} and \ref{s:wt4} where we list all known CM newforms of weight 3 and 4 with rational coefficients up
to twisting. This will answer our motivating question which
originated from singular K3 surfaces over $\Q$.

The proof of
Theorem \ref{1:1} is organised as follows: Denote the map of Theorem \ref{1:1}, which sends a
Hecke character to its CM-field, by $\jj$:
\begin{eqnarray*}
\jj:
\left\{\begin{matrix}
\text{Hecke characters}\\ \text{with $\infty$-type $l$ and rational coefficients}\\ \text{up to twisting}
       \end{matrix}\right\}
& \to & \left\{\begin{matrix}
\text{Imaginary quadratic fields $K$}\\ \text{with exponent}\;\, e_K\mid l
                                          \end{matrix}\right\}\\
\psi & \mapsto & K
\end{eqnarray*}
Then Theorem \ref{1:1} states that $\jj$ is well-defined (i.e.~all CM-fields occurring
have exponent $e_K\mid l$) and bijective. This will be established
as follows:

In the next section, we use an argument of Serre to check that
$\jj$ is well-defined (Proposition \ref{Prop:exp}). Then, section \ref{s:surj} gives the proof of the surjectivity of
$\jj$ (Corollary \ref{Cor:surj}). Using elementary genus theory,
we will also achieve a first step towards the injectivity of $\jj$
(Lemma \ref{Lem:minimal-twists}).
Finally, the proof of the injectivity of $\jj$ will be sketched in section
\ref{s:inj}.

\begin{Remark}\label{Rem:application}
Theorem \ref{1:1} can be used to determine a CM newform with rational coefficients
from very few explicit data. Given the bad primes, it suffices to check a small number of primes to deduce the
precise form. For instance this metod is used in \cite{Sch2}.
\end{Remark}

\section[Well-definedness]{$\boldsymbol{\jj}$ is well-defined}
\label{s:def}
\label{s:im}

In this section, we will prove that $\jj$ is well-defined.
Since $\jj$ is consistent with twisting, the following proposition provides an equivalent
formulation. It goes
back to an argument of Serre in \cite[Rem.~1.8]{L}:

\begin{Proposition}\label{Prop:exp}
Let $K$ be an imaginary quadratic field. Let $\psi$ be a
Hecke character of $K$ with rational coefficients and
$\infty$-type $l$. Then $e_K\mid l$.
\end{Proposition}

For the proof of the proposition, let $\m$ denote the conductor of $\psi$ and $G=(\OO_K/\m)^*$. We define the character
\begin{eqnarray*}
\eta:\;\, G & \to & \;\;\;\;\C[*]\\
\alpha & \mapsto & \frac{\psi(\alpha\,\OO_K)}{\alpha^l}.
\end{eqnarray*}
Formally, the definition
coincides with that of $\ez$, the character describing the
operation of $\psi$ on $\Z/M$. It follows from Corollary \ref{Cor:ez} that $\eta$ is determined on 
\[
G_{\Z} = \; G\cap \text{im} ((\Z/M)^* \to G) \; = \;\{\alpha\in G\,; \exists\; a\in\Z: \alpha\equiv a \mod{\m}\}.
\]
\begin{Corollary}\label{Cor:etaz}
If $\psi$ has real coefficients and $\infty$-type $l$, then \[
\eta|_{G_{\Z}}=\chi_K^l.\]
\end{Corollary}
Another essential ingredient for the proof of Proposition
\ref{Prop:exp} is provided by the next
\begin{Lemma}\label{Lem:image}
If the coefficients of $\psi$ are rational, then\[
\text{im } \eta \subseteq\OO_K^*.\]
\end{Lemma}
\emph{Proof:} Since $G$ is finite, $\eta$ has image some
roots of unity. In particular, these are algebraic integers. On
the other hand, $\eta$ has values in $K$ by Lemma \ref{Lem:rat}.
\qed

\begin{Corollary}\label{Cor:ideals}
Let $\psi$ be a Hecke character of $K$ with rational coefficients and
$\infty$-type $l$. Then, for any $\alpha\in\OO_K$ coprime to the
conductor $\m$,
\[
\psi(\alpha\,\OO_K)\,\OO_K=\alpha^l\,\OO_K=(\alpha\,\OO_K)^l
\]
as (principal) ideals of $\OO_K$.
\end{Corollary}

We can now prove Proposition
\ref{Prop:exp}: Let $\ia$ be an ideal of $\OO_K$ with order $n$
in $Cl(\OO_K)$, i.e.~$\ia^n=\alpha\,\OO_K$ for some $\alpha\in\OO_K$.
If $\ia$ is prime to $\m$, then Corollary \ref{Cor:ideals} gives
\[
(\psi(\ia)\,\OO_K)^n=(\alpha\,\OO_K)^l.
\]
Since both $\psi(\ia)\,\OO_K$ and $\alpha\,\OO_K$ are principal ideals, $n$
has to divide $l$ by the minimality of the order. This completes
the proof of Proposition \ref{Prop:exp}. \qed

The following corollary of Proposition \ref{Prop:exp} is a consequence of genus theory:
\begin{Corollary}[Odd $\boldsymbol\infty$-type]\label{Cor:odd}
Let $l\in\N$ be odd. If there is a Hecke character of $K$ with rational coefficients and $\infty$-type $l$, then
either $\Delta_K$ is prime (with $\Delta_K\equiv 3\mod 4$) or $\Delta_K=4$ or  $8$.
\end{Corollary}

\section[Surjectivity]{The surjectivity of $\boldsymbol{\jj}$}
\label{s:surj}

In this section, we prove the surjectivity of $\jj$ by constructing canonical Hecke characters for each imaginary quadratic field in question. As a first step towards the injectivity of $\jj$, we then show that these canonical Hecke characters are identified by twisting.

\begin{Lemma}\label{Lem:surj}
Let $K$ be an imaginary quadratic field and $l$ a
multiple of $e_K$. Then there is a Hecke character of $K$
with rational coefficients and $\infty$-type $l$. 
\end{Lemma}

\emph{Proof:} If $\Delta_K=3,4,8$ (or any other class number 1 discriminant), we can use the Hecke character associated to the elliptic curve with CM by $\OO_K$. Then consider its $l$-th power. 

In any other case, we construct Hecke characters with trivial conductor
if $l$ is even, resp.~with conductor $\m=(\sqrt{-\Delta})$ if $l$
is odd as follows: On the principal ideals of $\OO_K$, we define $\psi$ by
\[
\alpha\,\OO_K\mapsto (\chi_K(\text{Re}(\alpha)) \alpha)^l.
\]
For the non-principal ideals, we consider a factorisation of the
abelian group $Cl(\OO_K)$ into cyclic groups
\[
Cl(\OO_K)=C_{n_1} \times \hdots \times C_{n_r}.
\]
Then $\psi$ is completely determined by its operation on a set of
generators $\ia_1,\hdots,\ia_r$ (chosen coprime to
the conductor). Let $\mathfrak{a}_i^{n_i}=\alpha_i\,\OO_K$ for some
$\alpha_i\in\OO_K$. We obtain
\[
\psi(\mathfrak{a}_i)^{n_i}=\psi(\alpha_i\,\OO_K)=(\chi_K(\text{Re}(\alpha_i))\alpha_i)^l.
\]
In case $n_i$ is odd, this uniquely determines the root
\begin{eqnarray}\label{eq:odd}
\psi(\mathfrak{a}_i)=(\chi_K(\text{Re}(\alpha_i))\alpha_i)^{\frac{l}{n_i}}\in\OO_K.
\end{eqnarray}
If $n_i$ is even, there is a choice of sign in the image of $\mathfrak{a}_i$ under $\psi$:
\begin{eqnarray}\label{eq:even}
\psi(\mathfrak{a}_i)=\pm \alpha^{\frac{l}{n_i}}\in\OO_K.
\end{eqnarray}
Fix a sign for each generator $\ia_i$. This gives rise to a Hecke character with
rational coefficients and conductor as claimed.\qed

\begin{Corollary}\label{Cor:surj}
The map $\jj$ is surjective.
\end{Corollary}

\begin{Definition}
We call the Hecke characters constructed above \emph{canonical}.
\end{Definition}

\begin{Remark}\label{Rem:conductor}
It follows from the construction that we could equivalently have used the precise shape of the conductor to define the notion of a canonical Hecke character.
\end{Remark}

We conclude this section with a first step towards the
injectivity of $\jj$. 

\begin{Lemma}\label{Lem:minimal-twists}
The canonical Hecke characters of $K$ are equivalent under twisting.
\end{Lemma}

\emph{Proof:} If $e_K$ is odd, the canonical Hecke character is uniquely given by (\ref{eq:odd}). If $e_K$ is even, there are choices of sign due to (\ref{eq:even}). Changing signs is achieved by twisting with a fundamental character of $K$ (cf.~\cite[\S 3.B]{C}).\qed

\section[Injectivity]{The injectivity of $\boldsymbol{\jj}$}
\label{s:inj}
\label{s:unram-odd}

In this section, we will show how the injectivity of $\jj$ can be
proven. It will be derived from the following

\begin{Proposition}\label{Thm:twist}
Let $\psi$ be a Hecke character with rational coefficients. Then there is a twist of $\psi$ which is a canonical Hecke character.
\end{Proposition}

\begin{Corollary}\label{Cor:inj}
$\jj$ is injective.
\end{Corollary}

\emph{Proof of the Corollary:} Consider two Hecke characters of $K$ with rational coefficients and same $\infty$-type. By Prop.~\ref{Thm:twist}, there are twists which are canonical Hecke characters. By Lemma \ref{Lem:minimal-twists}, these twists are equivalent under twisting.\qed

By Remark \ref{Rem:conductor}, proving Prop.~\ref{Thm:twist} amounts to exhibiting a twist with the right conductor. Essentially we have to twist away the unramified primes. We will work with the character $\eta$ which was defined in section \ref{s:im}.
This suffices since twisting $\psi$ corresponds to
twisting $\eta$:
\begin{equation}\label{eq:twist}
(\psi, \eta) \leftrightarrow (\psi\otimes\Xi, \eta\otimes\Xi).
\end{equation}

The twisting in Prop.~\ref{Thm:twist} will be exhibited primewise. For every prime $\p$ dividing the conductor $\m$, we perform a twist such that the $\p$-parts of the conductor of the twist and of the canonical Hecke characters agree. We write $\p$ (and $\bar\p$) for the primes of $\OO_K$ and $p$ for the corresponding primes of $\Z$.

If a prime $\p$ divides $\m$, we factor $\m=\m'\p^{e_\p}$ or $\m=\m'\p^{e_\p}\bar\p^{e_{\bar\p}}$, such that $\m'$ is prime to $p$. Then we restrict $\eta$ to the multiplicative subgroup
\[
\Omega=\{\alpha\in(\OO_K/\m)^*: \alpha\equiv 1\mod{\m'}\}\subset G.
\]
We denote the restrictions
\[
\teta=\eta|_\Omega \;\;\; \text{ and }\;\;\; \Omega_{\Z}=\Omega\cap G_{\Z}.
\]
\begin{Lemma}\label{Lem:unram}
In the above notation,
\[
\teta|_{\Omega_{\Z}}=\begin{cases}
1, & \text{if $\p$ is unramified},\\
\chi_K^l, & \text{if $\p$ is ramified}..
                     \end{cases}
\]
\end{Lemma}
\emph{Proof:} If the $\infty$-type $l$ of $\psi$ is even or $\p$ is ramified, the statement follows from Corollary \ref{Cor:etaz}. If $l$ is odd, then we know that
\begin{equation}\label{here}
\eta|_{G_{\Z}}=\chi_K  \;\;\text{ (Cor.~\ref{Cor:etaz})}.
\end{equation}
Recall that in this setting $\Delta_K$ is either an odd prime or 4
or 8 (Cor.~\ref{Cor:odd}). In the first case, we can use the
property $(\sqrt{-\Delta_K})\mid\m$ from Corollary \ref{Cor:ez}.
If $\p$ is unramified, this gives $(\sqrt{-\Delta_K})\mid\m'$, so that
\[
a\in\Omega_{\Z} \Rightarrow a\equiv 1\mod{\Delta_K}.
\]
Since $\teta|_{\Omega_{\Z}}=(\eta|_{G_{\Z}})|_{\Omega_{\Z}}$, equation (\ref{here}) now gives the claim. For $\Delta_K=4$ or $8$, the same argument can be applied using Remark \ref{Rem:ez}. \qed

By the Chinese remainder theorem,
\[
\Omega\cong(\OO_K/\p^{e_\p})^* \;\;\text{ or } \;\;\Omega\cong(\OO_K/\p^{e_\p})^* \oplus (\OO_K/\bar\p^{e_{\bar\p}})^*.
\]
The further factorisation of this abelian group will allow us to determine $\teta$ explicitly. Then we will exhibit the corresponding twist which produces the right conductor. This will also give us the possible values for the exponents $e_\p$. We collect them in Proposition \ref{Prop:bounds-even}.

Here we will only consider the case where $p>2$ and $\OO_K^*=\{\pm 1\}$. The other cases, $p=2$ or $\Delta_K=3$ or $4$, require slightly more attention since there are additional characters to be considered. However, no essential new ideas are needed. 
The details can be found in \cite{S-Diss}. 

To achieve Prop.~\ref{Thm:twist} in the above case, we apply the following lemma which makes the twisting explicit (cf.~(\ref{eq:twist})):

\begin{Lemma}\label{Prop:unram-odd}\label{Lem:ram}
Let $p>2$ and $\Delta_K\neq 3, 4$. If $p$ does not ramify in $K$, then $\teta=1$ or $\chi_K\circ\NN$. If $p$ ramifies, then $\teta$ is uniquely given by Lemma \ref{Lem:unram}.
\end{Lemma}

If $(p,\m)=1$, then $\teta=1$, so there is nothing to prove. If $(p, \m)\neq 1$, we distinguish between the splitting behaviour of $p$ in $K$ to prove Lemma \ref{Lem:ram}.

\subsection[$p$ is inert in $K$]{$\boldsymbol{p}$ is inert in $\boldsymbol K$}
\label{s:inert}

In the inert case, $\teta$ operates on the quotient $\Omega\cong(\OO_K/p^{e})^*$ where $e=e_p$. As an abelian group, it has a factorisation
\[
\Omega\cong(\OO_K/p^{e})^* \cong C_{p^2-1} \times \tilde\Omega
\]
where $\tilde\Omega$ is a group of cardinality $p^{2(e-1)}$. In particular, this cardinality is odd. Hence $\teta$ as a quadratic character operates trivially on $\tilde\Omega$. This gives $e\leq 1$.

We deduce that $\Omega=(\OO_K/p)^*=C_{p^2-1}$ with a non-trivial action of the quadratic character $\teta$. Let $\chi_p$ denote the unique quadratic Dirichlet character of conductor $p$. Then $\chi_p\circ\NN$ defines another non-trivial quadratic character on the cyclic group $\Omega$. Hence the two characters $\teta$ and $\chi_p\circ\NN$ coincide.
\qed

\subsection[$p$ splits in $K$]{$\boldsymbol{p}$ splits in $\boldsymbol K$}
\label{s:unram-split}

In the split case, we write $p=\p\bar\p$ such that
$\m=\m'\p^{e_\p}\bar\p^{e_{\bar\p}}$. Hence 
\[
\Omega\cong (\OO_K/\p^{e_\p})^* \oplus (\OO_K/\bar\p^{e_{\bar\p}})^* \cong \Omega_1\oplus\Omega_2.
\]
The respective quotient maps $\Omega\to\Omega_i$ shall be denoted by $[\cdot]_i$. Since $\p\neq\bar\p$, the summands are isomorphic to $(\Z/p^{e_\p})^*$ and $(\Z/p^{e_{\bar\p}})^*$, respectively.

We write $\tilde{\eta}=\teta_1\cdot \teta_2$ with the characters $\teta_i$ operating on $\Omega_i$. We can view them as acting on $(\Z/p^e)^*$ for the respective exponent $e$. With this notion, our next result is that the characters are conjugate:
\begin{eqnarray}\label{eq:conj}
\teta_1=\overline{\teta_2}.
\end{eqnarray}
This is a direct consequence of the equality from Lemma \ref{Lem:unram}
\[
1=\teta(a)=\teta_1([a]_1)\, \teta_2([a]_2)\;\;\forall\; a\in\Omega_{\Z}.
\]
To deduce the result, we only need that $[a]_1\equiv[a]_2\mod{p^{\text{min}(e_\p,e_{\bar\p})}}$ for all $a\in\Omega_{\Z}$. In particular, the conjugacy (\ref{eq:conj}) implies that
\begin{eqnarray}\label{eq:exp-odd}
e=e_\p=e_{\bar\p}.
\end{eqnarray}
Our next claim is that the characters $\teta_i$ are in fact
quadratic; then, by conjugation
\[
\teta_1=\teta_2.
\]
To prove this, let us assume on the contrary that the character $\tilde{\eta}_1$ is not quadratic. We shall establish a contradiction.
By assumption, there is an element $\alpha\in\Omega$ with $\tilde{\eta}_1([\alpha]_1)\not\in\{\pm 1\}$. Let $n_{\p}, n_{\m'}$ be the respective orders of $\p$ and $\m'$ in the class group $Cl(K)$, so that $\p^{n_{\p}}=(\pi)$ and $\m'^{\,n_{\m'}}=(\mu)$  
for some $\pi, \mu\in\OO_K$. We then consider the elements
\[
\alpha+k\pi\mu\in\Omega \;\;\; \text{ for }\; k\in\Z.
\]
Taking $k$ in the range of $0,\hdots,p^e-1$, the residue classes
$[\alpha+k\pi\mu]_2$ run through the whole of $\Z/p^e$ by definition. Hence there is a $k_0$ such that 
$[\alpha+k_0\pi\mu]_2=[1]_2$. This gives the required contradiction, since $\teta$ is quadratic by Lemma \ref{Lem:image}:
\begin{eqnarray*}
\tilde{\eta}(\alpha+k_0\pi\mu) & = & \tilde{\eta}_1([\alpha+k_0\pi\mu]_1)\; 
\tilde{\eta}_2([\alpha+k_0\pi\mu]_2)\\
& = & \tilde{\eta}_1([\alpha]_1)\; \tilde{\eta}_2([1]_2)=\tilde{\eta}_1([\alpha]_1)\not\in\{\pm 1\}.
\end{eqnarray*}
Thus  the $\teta_i$ are quadratic. This implies $e\leq 1$, since
\[
\Omega_i\cong(\Z/p^e)^*\cong C_{p-1}\times C_{p^{e-1}}
\]
with the second factor of odd cardinality. In particular, $\teta_i$ is a non-trivial quadratic character on $\F_p^*$. Hence $\teta_i$ $(i=1,2)$ coincides with $\chi_p$, the unique quadratic Dirichlet character of conductor $p$. 

This property implies the claim $\teta=\chi_p\circ\NN$ of Lemma \ref{Lem:ram}. To see this, we note that $[\alpha]_2=[\bar\alpha]_1$ in $\F_p^*$;
in particular, for $\alpha\in\Omega_{\Z}$, we can omit the
subscript, regardless of the factor. For general $\alpha\in\Omega$
we thereby obtain
\begin{eqnarray*}
\tilde{\eta}(\alpha) & = & \tilde{\eta_1}([\alpha]_1)\tilde{\eta_2}([\alpha]_2)\\
& = & \chi_p([\alpha]_1) \chi_p([\bar{\alpha}]_1)\\
& = & \chi_p([\NN(\alpha)])\\
& = & \chi_p(\NN(\alpha)) = \chi_p\circ\NN(\alpha).
\end{eqnarray*}
This completes the proof of Lemma \ref{Lem:ram} in the split case. \qed

\subsection[$p$ ramifies in $K$]{$\boldsymbol{p}$ ramifies in $\boldsymbol{K}$}
\label{s:ram-odd}

In the ramified case, $\teta$ is a quadratic character on
\[
\Omega \cong (\OO_K/\p^e)^* \cong C_{p-1} \times \tilde\Omega
\]
where $\tilde\Omega$ denotes a group of odd cardinality $p^{e-1}$.
Hence again $e\leq 1$. But then,
\[
\Omega \cong (\OO_K/\p)^*\cong (\Z/p)^* \cong \Omega_{\Z}.
\]
In consequence, $\teta$ is uniquely determined by Lemma \ref{Lem:unram}. This completes the proof of Lemma \ref{Lem:ram}. Prop.~\ref{Thm:twist} and thereby the injectivity of $\jj$ and Thm.~\ref{1:1} follow. \qed

\section{Exponents}
\label{s:bounds}

We have seen that any Hecke character with rational coefficients admits a twist which has the canonical conductor from Lemma \ref{Lem:surj}. Since all canonical Hecke characters are equivalent under twisting, this proves Thm.~\ref{1:1}. In this section, we collect the possible exponents of the conductor. We also put them in relation with results of Serre. 

\begin{Proposition}\label{Prop:bounds-even}
Let $K$ be an imaginary quadratic field. Let $\psi$ be a Hecke character of $K$ with rational coefficients and $\infty$-type $l$. Let $\p$ be a prime of $\OO_K$. Denote the exponent of the conductor of $\psi$ at $\p$ by $e_\p$. Then $e_\p=e_{\bar\p}$ with the following possible values:

\begin{center}
\begin{tabular}{|c||c|c|c||c|c|c|c|c|}
\hline
prime  & \multicolumn{3}{c||}{$\p$ prime to 2} & \multicolumn{5}{c|}{$\p$ above 2\;\;\;} \\
\hline
\hline ramification & \multicolumn{2}{c|}{ramified\;\;} & unram. & \multicolumn{2}{c|}{$4||\Delta_K$} & \multicolumn{2}{c|}{$8|\Delta_K$} & unram. \\
\hline $\infty$-type & odd & even & arb. & odd & even & odd & even & arb.\\
\hline
\hline $\Delta_K\neq 3,4$ & 1 & 0 & 0,1  & - & 0,4 & 5 & 0,2 & 0,2,3  \\
\hline $\Delta_K=4$ & - & - & 0,1 & 3,4,6 & 0,2,4,6 & - & - & -  \\
\hline $\Delta_K=3$ & 1,2,4 & 0,2,4 & 0,1  & - & - & - & - & 0,1,2,3 \\
\hline
\end{tabular}
\end{center}
\end{Proposition}

Thm.~\ref{Thm:Shimura} gives the corresponding statement for the level of the associated newform. The possible levels can be compared to results of Serre \cite[(4.8.8) \& Application \`a (4.8.8) b)]{Se}. They agree with the maximal exponents for the conductor $N=\prod p^{e_p}$ of a modular integral
two-dimensional Galois representation of the middle cohomology group of an odd-dimensional
smooth projective variety over $\Q$:
\[
\begin{matrix}
e_p\leq 2, & \text{ if } p>3,\\
e_p\leq 5, & \text{ if } p=3,\\
e_p\leq 8, & \text{ if } p=2.
\end{matrix}
\]
Note that Serre's bounds only apply to even weight, while our results only hold for CM-forms, but of arbitrary weight.

\section{The CM newforms of weight 3}
\label{s:wt3}

In this section we will give all CM newforms of weight 3 with
rational Fourier coefficients up to twisting (except for possibly
one). By a theorem of Livn\'e \cite{L}, these constitute exactly
the modular forms which can a priori appear as non-trivial factor
$L(T_X,s)$ in the $L$-series of singular K3 surfaces over $\Q$.
Hence this section answers the question which was our original
motivation. In particular, we deduce the following result:

\pagebreak

\begin{Proposition}\label{Lem:Class-d}
Consider the following classifications of singular K3 surfaces over $\Q$:
\begin{enumerate}[(i)]
\item\label{(i)} by the discriminant $d$ of the transcendental
lattice of the surface up to square, 
\item\label{NS}
by the discriminant $-d$ of the N\'eron-Severi lattice of the surface up to square,
\item\label{(ii)} by the
associated newform up to twisting, 
\item\label{(iii)} by the level
of the associated newform up to square, 
\item\label{(iv)} by the
CM-field $\Q(\sqrt{-d})$ of the associated newform.
\end{enumerate}
These classifications are equivalent. In particular, $\Q(\sqrt{-d})$ has exponent 1 or 2.
\end{Proposition}

\emph{Proof:} 
The relation between the discrimininants of transcendental lattice in (\ref{(i)}) 
and N\'eron-Severi lattice in (\ref{NS}) is due to Nikulin \cite{N}.
The classifications by (\ref{(ii)}) and (\ref{(iv)})
coincide due to Theorem \ref{1:1} which also gives the extra
claim. Hence, the equivalence of
(\ref{(ii)}) and (\ref{(iv)}) with (\ref{(i)}) follows from
\cite[Thm.~1.3, Ex.~1.6]{L}. Finally, the equivalence
(\ref{(iii)})$\Leftrightarrow$(\ref{(iv)}) can be read off from
Proposition \ref{Prop:bounds-even}. \qed

\begin{Remark}
The Torelli theorem \cite{PSS} identifies singular K3 surfaces with their transcendental lattices up to isomorphism.
The intersection form of the transcendental lattice is an even positive definite binary quadratic form. This provides an explicit link to class group theory.

The conclusion of class group exponent 1 or 2 for singular K3 surfaces over $\Q$ in Prop.~\ref{Lem:Class-d} has been established independently in \cite{S-fields}. The reasoning in \cite{S-fields} is mainly geometrical and builds on work by Shioda--Inose on singular K3 surfaces \cite{SI} and by Shioda--Mitani on singular abelian surfaces \cite{SM}.
\end{Remark}

For newforms of weight 3 with rational coefficients, we have to consider imaginary quadratic fields whose class group consists only of 2-torsion. These fields are related to Euler's \emph{idoneal numbers}. They  were also studied by Gauss. Weinberger showed that the list of 65 known fields
is almost complete \cite[Thm.~2]{Wb}. More precisely, he stated that there is at
most one imaginary quadratic field with exponent 2 and $\Delta_K>5460$. He showed that this field would have
$\Delta_K>2\cdot 10^{11}$. Subject to a Siegel-Landau condition on
$L$-functions for odd real Dirichlet characters, which in particular would be a consequence of ERH  (cf.~\cite[Thm.~4]{Wb}),
he proved that the known fields with $\Delta_K\leq 5460$ constitute the
complete list.

In this section, we will give a newform $f=\sum_n a_nq^n$ of weight 3 with rational Fourier coefficients and CM by $K$ for every known field $K$. The Fourier coefficients  were calculated by a computer program after fixing a choice of sign in (\ref{eq:even}). In Table \ref{T:wt3}, the first column lists the level $N$ of the
respective CM newforms. By Proposition \ref{Lem:Class-d}, such a
newform has CM by $\Q(\sqrt{-N})$. This also gives its nebentypus.
We then list the Fourier coefficients at the first 30 primes.

\begin{table}[ht!]

\begin{scriptsize}

\begin{tabular}{|c||c|c|c|c|c|c|c|c|c|c|c|c|c|c|c|c|c|c|}
\hline
$N$ & 2 & 3 & 5 & 7 & 11 & 13 & 17 & 19 & 23 & 29 & 31 & 37 & 41 & 43 & 47 & 53\\
\hline
\hline

7&-3&0&0&-7&-6&0&0&0&18&-54&0&-38&0&58&0&-6
\\ \hline
8&-2&-2&0&0&14&0&2&-34&0&0&0&0&-46&14&0&0
\\ \hline
11&0&-5&-1&0&-11&0&0&0&35&0&-37&-25&0&0&50&-70
\\ \hline
15&-1&3&-5&0&0&0&14&-22&-34&0&2&0&0&0&14&86
\\ \hline
16&0&0&-6&0&0&10&-30&0&0&42&0&-70&18&0&0&90
\\ \hline \hline
19&0&0&-9&-5&3&0&15&-19&-30&0&0&0&0&-85&75&0
\\ \hline
20&2&-4&-5&4&0&0&0&0&-44&-22&0&0&62&76&4&0
\\ \hline
24&2&-3&-2&-10&10&0&0&0&0&-50&38&0&0&0&0&94
\\ \hline
27&0& 0& 0& -13& 0& -1& 0& 11& 0& 0& -46& 47& 0& -22& 0& 0
\\ \hline
35&0&-1&5&-7&-13&19&-29&0&0&23&0&0&0&0&31&0
\\ \hline \hline
40&2&0&-5&-6&-18&6&0&-2&26&0&0&54&-78&0&-86&-74
\\ \hline
43&0&0&0&0&-21&-17&-9&0&3&0&19&0&39&-43&-78&63
\\ \hline
51&0&3&-7&0&5&-25&-17&-13&29&-10&0&0&65&35&0&0
\\ \hline
52&2&0&0&-12&-4&-13&-18&12&0&6&36&0&0&0&68&-102
\\ \hline
67&0&0&0&0&0&0&-33&-29&-21&-9&0&7&0&0&27&0
\\ \hline \hline
84&2&3&-4&-7&-20&0&20&10&4&0&-50&-10&68&0&0&0
\\ \hline
88&2&0&0&0&-11&-18&0&-6&-42&14&-26&0&0&42&6&0
\\ \hline
91&0&0&3&-7&0&13&0&-25&-45&-33&55&0&-30&-5&-81&15
\\ \hline
115&0&0&5&-9&0&0&11&0&-23&-57&-53&51&-33&-6&0&-101
\\ \hline
120&2&3&-5&0&2&-14&-26&0&-14&38&-58&34&0&-74&34&0
\\ \hline\hline
123&0&3&0&0&-19&0&-7&0&0&17&-61&-49&-41&-37&53&-58
\\ \hline
132&2&3&0&-8&-11&0&-32&16&2&-8&0&-58&16&64&-82&0
\\ \hline
148&2&0&0&0&0&0&0&-36&-28&0&-12&-37&-66&12&0&-42
\\ \hline
163&0&0&0&0&0&0&0&0&0&0&0&0&-81&-77&-69&-57
\\ \hline
168&2&3&0&-7&0&-2&-22&0&-38&-26&34&0&26&-82&0&22
\\ \hline\hline
187&0&0&0&-3&11&0&-17&0&0&41&0&0&-71&0&-93&-81
\\ \hline
195&0&3&5&1&-17&-13&-31&0&-19&0&0&61&43&0&0&41
\\ \hline
228&2&3&0&0&-16&0&0&-19&8&-56&-14&0&-32&0&56&-8
\\ \hline
232&2&0&0&0&0&0&0&0&0&-29&-54&-42&0&0&-22&0
\\ \hline
235&0&0&5&0&0&-21&0&0&-1&0&0&0&0&39&-47&0
\\ \hline\hline
267&0&3&0&0&0&0&0&0&-43&-31&0&0&-7&0&0&0
\\ \hline
280&2&0&5&-7&0&0&6&-18&0&0&0&-66&0&-54&66&-34
\\ \hline
312&2&3&0&0&0&-13&0&-14&0&-46&0&22&-74&0&-62&2
\\ \hline
340&2&0&5&0&-12&0&-17&0&0&0&28&6&0&-84&-76&0
\\ \hline
372&2&3&0&0&0&0&-28&0&0&-4&-31&0&0&-38&-92&44
\\ \hline\hline
403&0&0&0&0&-9&13&0&0&0&0&-31&43&0&0&0&0
\\ \hline
408&2&3&0&0&0&0&-17&0&-22&0&0&-62&14&0&0&-98
\\ \hline
420&2&3&5&-7&8&-16&0&-32&0&0&-8&0&-2&-26&-46&-104
\\ \hline
427&0&0&0&7&0&0&-27&0&0&0&1&0&0&0&0&0
\\ \hline
435&0&3&5&0&-7&0&0&0&-41&-29&0&-71&53&-59&0&19
\\ \hline\hline
483&0&3&0&7&-1&0&0&-31&-23&0&0&0&-79&0&-67&83
\\ \hline
520&2&0&5&0&0&-13&0&0&-6&0&-42&0&0&0&0&54
\\ \hline
532&2&0&0&7&0&-12&0&-19&0&0&0&0&44&0&18&0
\\ \hline
555&0&3&5&0&0&-11&0&0&0&-53&0&-37&0&49&-91&-79
\\ \hline
595&0&0&5&7&0&9&-17&0&-39&0&-57&-11&-37&0&-59&0
\\ \hline\hline
627&0&3&0&0&11&7&-23&-19&0&0&0&0&0&0&0&-103
\\ \hline
660&2&3&5&0&-11&4&0&-28&0&-52&0&0&-28&0&0&-26
\\ \hline
708&2&3&0&0&0&0&0&0&0&0&-56&0&0&-32&0&0
\\ \hline
715&0&0&5&0&11&-13&21&-27&0&0&0&-69&17&-31&-49&0
\\ \hline
760&2&0&5&0&0&0&0&-19&0&-18&0&0&0&-66&0&0
\\ \hline\hline
795&0&3&5&0&0&0&-19&0&0&0&0&0&-77&0&41&-53
\\ \hline
840&2&3&5&-7&0&0&0&0&0&2&-22&46&-58&0&-74&0
\\ \hline
1012&2&0&0&0&11&0&-12&0&-23&0&0&0&0&0&0&0
\\ \hline
1092&2&3&0&7&0&-13&8&0&-32&0&0&0&0&0&-88&0
\\ \hline
1155&0&3&5&7&-11&0&1&-17&-31&47&0&0&0&-79&0&29
\\ \hline\hline
1320&2&3&5&0&-11&0&0&0&0&0&0&-14&38&-46&0&0
\\ \hline
1380&2&3&5&0&0&0&0&-8&-23&0&0&-64&0&0&0&0
\\ \hline
1428&2&3&0&7&0&0&-17&4&0&-44&0&0&0&0&0&0
\\ \hline
1435&0&0&5&7&0&0&0&-3&0&0&0&0&-41&0&0&-99
\\ \hline
1540&2&0&5&7&-11&0&0&0&24&0&-48&0&-72&0&0&0
\\ \hline\hline
1848&2&3&0&7&-11&0&0&0&0&0&0&0&0&-2&-38&62
\\ \hline
1995&0&3&5&7&0&0&0&-19&-11&-37&43&-59&0&0&0&0
\\ \hline
3003&0&3&0&7&11&-13&0&0&0&19&-29&0&-61&0&0&0
\\ \hline
3315&0&3&5&0&0&13&-17&0&0&7&-23&0&0&0&0&89
\\ \hline
5460&2&3&5&7&0&-13&0&0&0&0&0&-4&0&-44&0&-76\\
\hline
\end{tabular}
\vspace{0.1cm}
\caption{The CM newforms of weight 3 with rational Fourier coefficients}
\label{T:wt3}
\end{scriptsize}
\end{table}

\begin{table}
\begin{scriptsize}

\begin{tabular}{|c|c|c|c|c|c|c|c|c|c|c|c|c|c|}
\hline
59 & 61 & 67 & 71 & 73 & 79 & 83 & 89 & 97 & 101 & 103 & 107 & 109 & 113\\
\hline
\hline
0&0&-118&114&0&-94&0&0&0&0&0&186&106&-222
\\ \hline
-82&0&62&0&-142&0&158&146&-94&0&0&-178&0&98
\\ \hline
107&0&35&-133&0&0&0&-97&95&0&-190&0&0&215
\\ \hline
0&-118&0&0&0&98&-154&0&0&0&0&-106&-22&206
\\ \hline
0&-22&0&0&-110&0&0&-78&130&-198&0&0&-182&-30
\\ \hline\hline
0&103&0&0&-25&0&90&0&0&-102&0&0&0&0
\\ \hline
0&-58&-116&0&0&0&76&-142&0&122&-44&124&38&0
\\ \hline
10&0&0&0&50&-58&-134&0&-190&190&-10&-86&0&0
\\ \hline
0& -121&
-109& 0& -97& 131& 0& 0& 167& 0& -37& 0& -214& 0
\\ \hline
0&0&0&2&34&-157&-86&0&-149&0&199&0&-97&0
\\ \hline\hline
78&0&0&0&0&0&0&18&0&0&186&0&0&0
\\ \hline
-54&0&91&0&0&-14&123&0&-193&159&-181&42&-169&0
\\ \hline
0&0&-70&-130&0&0&0&0&0&0&155&-211&0&-199
\\ \hline
-116&-86&108&-92&0&0&-68&0&0&-6&0&0&0&174
\\ \hline
51&0&-67&-126&79&0&-102&111&0&0&-62&147&0&0
\\ \hline\hline
0&0&0&100&0&0&0&-172&0&-148&94&-164&-118&0
\\ \hline
0&78&0&54&0&0&122&-174&-158&-194&118&-182&174&-126
\\ \hline
90&0&0&0&-29&67&159&-165&131&0&0&-150&0&135
\\ \hline
3&0&111&27&0&0&-41&0&-174&87&114&191&0&19
\\ \hline
98&0&-26&0&0&38&0&0&0&182&0&0&0&166
\\ \hline\hline
0&-1&0&101&23&0&0&14&0&161&83&0&0&0
\\ \hline
74&0&0&-34&0&136&0&0&62&136&0&0&0&0
\\ \hline
44&0&0&0&-2&84&0&0&0&54&132&0&0&0
\\ \hline
0&-41&0&-21&0&0&3&0&31&0&0&0&0&63
\\ \hline
-106&94&-34&58&0&0&-58&122&0&0&178&0&0&0
\\ \hline\hline
-69&54&-53&0&129&-114&0&-9&0&0&19&61&-207&0
\\ \hline
-38&-73&-74&103&94&-37&0&-173&181&0&0&149&0&-34
\\ \hline
0&-106&58&0&-82&-146&128&64&0&0&-98&0&0&112
\\ \hline
-114&6&-98&0&0&42&-66&0&0&86&0&-18&0&0
\\ \hline
-117&-113&-54&-93&99&-77&0&-57&0&-33&0&-209&0&179
\\ \hline\hline
29&0&-133&0&-121&-109&77&-89&-73&-154&0&0&-49&137
\\ \hline
62&-102&-6&-138&-106&-122&0&0&166&-22&-46&74&0&0
\\ \hline
0&0&82&-14&0&-154&0&22&0&98&-106&-202&166&0
\\ \hline
0&0&-36&108&-126&-148&-4&-162&126&-138&36&0&0&-46
\\ \hline
-68&0&0&-44&0&34&0&116&-178&0&0&28&-154&0
\\ \hline\hline
0&0&0&0&-133&0&42&147&0&-201&-197&-189&0&-177
\\ \hline
-86&-14&0&74&0&0&-38&0&0&-2&-202&0&82&158
\\ \hline
0&0&106&128&104&0&26&-158&-184&118&0&0&-202&16
\\ \hline
57&-61&0&0&0&0&0&-66&0&141&0&-213&-209&-201
\\ \hline
0&0&0&0&1&0&79&62&49&173&0&-134&-217&0
\\ \hline\hline
-43&53&0&0&0&0&0&0&-82&41&137&122&0&-142
\\ \hline
0&0&-126&38&-114&0&-94&0&-66&0&154&0&-198&0
\\ \hline
0&0&-132&-124&0&-108&-138&-164&156&0&0&-52&0&0
\\ \hline
7&0&0&0&0&0&-19&67&46&0&169&29&0&0
\\ \hline
0&3&0&0&0&0&149&0&0&0&-66&129&0&-114
\\ \hline\hline
-91&0&0&-67&0&139&109&-31&0&-26&0&0&-86&17
\\ \hline
-58&0&-86&98&124&92&-164&0&0&92&-14&-116&0&94
\\ \hline
-59&0&16&-94&0&0&0&-176&0&-152&88&-22&0&-128
\\ \hline
0&0&-9&0&0&0&0&0&51&0&0&6&-42&0
\\ \hline
42&0&0&0&0&0&14&0&-186&0&-174&0&142&-154
\\ \hline\hline
0&0&-131&-17&-119&0&0&0&0&43&-59&161&0&14
\\ \hline
-22&38&0&-82&-134&0&0&38&-86&0&0&-206&0&-194
\\ \hline
0&76&42&0&0&0&0&0&0&0&-162&0&-196&0
\\ \hline
-64&0&-22&0&-62&0&-16&0&142&176&-158&136&0&0
\\ \hline
-113&67&0&0&0&0&-131&-53&-191&0&-179&0&0&149
\\ \hline\hline
0&-98&0&-122&14&-62&0&0&0&0&-146&0&-2&-214
\\ \hline
-112&0&0&-88&0&112&74&-98&56&0&0&-154&0&0
\\ \hline
0&-116&0&0&-92&-46&0&-94&-44&134&172&0&0&124
\\ \hline
0&0&-71&0&-141&0&-121&137&0&-167&-81&0&0&0
\\ \hline
8&-32&-64&0&0&-18&-54&0&-114&48&0&0&0&0
\\ \hline\hline
0&0&0&0&-118&0&-142&46&0&-106&0&0&-134&0
\\ \hline
0&0&1&47&-139&0&0&0&0&-197&0&0&0&0
\\ \hline
0&109&0&-131&0&0&23&0&103&0&0&-137&-211&0
\\ \hline
-103&0&-121&0&0&0&0&-43&0&0&0&0&133&0
\\ \hline
92&0&0&0&0&0&0&0&0&-188&0&0&0&44\\
\hline
\end{tabular}
\vspace{0.5cm}

\begin{center}
\begin{normalsize}
{\sf Tab. 1:} The CM newforms of weight 3 with rational Fourier
coefficients
\end{normalsize}
\end{center}

\end{scriptsize}

\end{table}

\section{The CM newforms of weight 4}
\label{s:wt4}

This section is concerned with CM newforms of weight 4 with
rational Fourier coefficients. By Theorem
\ref{1:1} we consider
imaginary quadratic fields whose class group is at most
3-torsion. The number of these fields is known to be finite. This
goes back to Boyd and Kisilevski \cite{BK}, but can be also be
found in Weinberger's article \cite{Wb}. Since we have not found an explicit list, we searched for them with
a straight forward computer program. Among the imaginary quadratic
fields with $\Delta_K\leq 100000$, we found 26 such
fields. 

For each such $K$, Table \ref{T:wt4} lists a
newform with CM by $K$, rational Fourier coefficients and minimal
level $N$. This newform is unique unless $K=\Q(\sqrt{-2})$ where the twist by $\chi_{-1}$ has
the same level. In the table, the level constitutes the
first column and is expressed as $N=M\cdot\Delta_K^2$. The Fourier
coefficients at the first 25 primes follow.

\pagebreak
\pagestyle{empty}

\begin{sidewaystable}
\centering
\begin{scriptsize}
\begin{tabular}{|c||c|c|c|c|c|c|c|c|c|c|c|c|c|c|c|c|c|c|c|c|c|c|c|c|c|}
\hline
$N$ & 2 & 3 & 5 & 7 & 11 & 13 & 17 & 19 & 23 & 29 & 31 & 37 & 41 & 43 & 47 & 53 & 59 & 61 & 67 & 71 & 73 & 79 & 83 & 89 & 97\\
\hline
\hline
$3^2$&0&0&0&20&0&-70&0&56&0&0&308&110&0&-520&0&0&0&182&-880&0&1190&884&0&0&-1330\\ \hline
$2\cdot 4^2$&0&0&22&0&0&-18&-94&0&0&-130&0&214&-230&0&0&518&0&830&0&0&1098&0&0&-1670&594\\ \hline
$7^2$&-5&0&0&0&-68&0&0&0&-40&-166&0&450&0&-180&0&590&0&0&-740&688&0&-1384&0&0&0\\ \hline
$4\cdot 8^2$&0&10&0&0&-18&0&90&106&0&0&0&0&-522&-290&0&0&846&0&-70&0&430&0&-1350&-1026&-1910\\ \hline
$11^2$&0&8&18&0&0&0&0&0&-108&0&340&-434&0&0&-36&-738&-720&0&-416&612&0&0&0&1674&-34\\ \hline\hline
$19^2$&0&0&14&-36&40&0&14&0&212&0&0&0&0&128&364&0&0&630&0&0&1078&0&-112&0&0\\ \hline
$23^2$&3&4&0&0&0&-74&0&0&0&282&-344&0&426&0&48&0&-396&0&0&1176&-1226&0&0&0&0\\ \hline
$31^2$&1&0&2&16&0&0&0&156&0&0&0&0&-278&0&-616&0&-740&0&684&1000&0&0&0&0&-1906\\ \hline
$43^2$&0&0&0&0&32&-90&-130&0&-140&0&108&0&22&0&500&-130&-904&0&-360&0&0&-1116&-680&0&290\\ \hline
$59^2$&0&-7&-21&-29&0&0&126&-119&0&-159&0&0&-525&0&0&327&0&0&0&-180&0&-835&0&0&0\\ \hline\hline
$67^2$&0&0&0&0&0&0&50&-144&-220&-266&0&270&0&0&-220&0&-104&0&0&788&-90&0&-1480&-374&0\\ \hline
$83^2$&0&5&0&-25&57&0&-75&0&180&309&-317&-335&378&0&0&0&471&-943&0&0&0&0&0&0&0\\ \hline
$107^2$&0&-1&0&0&-9&-11&0&-137&189&234&0&353&279&0&540&-711&0&911&0&0&0&-1127&-360&-1467&0\\ \hline
$139^2$&0&0&19&-11&43&93&0&0&0&-155&-345&-306&490&0&644&0&0&0&-29&575&0&-1091&1433&-1549&0\\ \hline
$163^2$&0&0&0&0&0&0&0&0&0&0&0&0&122&-360&-580&-770&0&918&0&-1012&0&0&-1040&0&990\\ \hline\hline
$211^2$&0&0&17&0&7&-83&0&165&0&0&0&-331&0&507&461&158&760&0&0&425&1242&-1285&1408&0&0\\ \hline
$283^2$&0&0&0&33&-71&79&0&0&-29&223&0&0&-61&0&0&0&-163&691&0&212&-630&0&1400&1423&537\\ \hline
$307^2$&0&0&0&12&-64&0&2&128&0&0&0&398&-506&0&0&34&0&0&0&-1196&0&-684&-1120&1526&1890\\ \hline
$331^2$&0&0&13&0&0&0&-139&-25&0&0&333&0&0&277&0&-707&0&0&-639&-1013&0&1141&248&1210&0\\ \hline
$379^2$&0&0&11&0&0&0&0&155&53&0&0&429&-197&0&0&0&0&-907&889&0&0&609&-1403&0&-846\\ \hline\hline
$499^2$&0&0&-1&0&0&0&0&0&0&115&83&0&0&313&389&0&0&0&-1029&0&1173&0&0&0&0\\ \hline
$547^2$&0&0&0&0&-56&-6&0&88&0&298&0&0&0&0&-388&766&0&0&-488&0&1114&0&0&0&-1758\\ \hline
$643^2$&0&0&0&27&0&0&0&0&-131&-257&259&0&0&0&0&751&0&0&0&0&0&0&931&-1183&-1897\\ \hline
$883^2$&0&0&0&0&0&-29&-137&0&0&-233&-111&0&0&0&0&0&0&0&1047&-1151&0&147&809&0&0\\ \hline
$907^2$&0&0&0&0&0&-25&0&69&65&0&0&0&-503&0&0&-665&0&0&0&1187&0&0&0&1501&0\\ \hline
$4027^2$&0&0&0&0&0&-69&125&-153&0&-182&0&0&134&416&0&-329&0&-43&-1095&1111&1109&443&-1511&-1649&1842\\ \hline
\end{tabular}

\caption{The CM newforms of weight 4 with rational Fourier coefficients}
\label{T:wt4}

\end{scriptsize}
\end{sidewaystable}



\section{Geometric Realisations}

Originally, we wanted to understand the modularity of singular K3
surfaces in terms of the associated newforms of weight 3. This has
been addressed by Prop.~\ref{Lem:Class-d} and Table \ref{T:wt3} in section \ref{s:wt3}. The
natural question now is which of these newforms actually occur in
the $L$-series of some K3 surface over $\Q$. We formulate this problem in general form:

\begin{Question}\label{Question}
Which CM newforms with rational coefficients (and weight $k$) have
geometric realisations in a smooth projective variety $X$ over $\Q$
with $h^{k-1,0}(X)=\dim H^0(X,\Omega_X^{k-1})=1$?
\end{Question}

In other words we ask whether the Galois representation associated
to a newform occurs as 2-dimensional motive of some smooth projective
variety which is defined over $\Q$. Note that the classical $(k-2)$-fold fibre product
of the universal elliptic curve with level $N$-structure which was
fundamental to Deligne's construction of the associated Galois
representation for a newform of weight $k$ and level $N$, in
general does not meet the given criterion.

What makes our situation special, is that we have very precise
associated geometric objects, different from Deligne's fibre
products. Namely, these objects are the elliptic curves with CM by
the respective field $K$. By the classical theory of CM, this
gives a positive answer to our question for weight $k=2$. For
higher weight, it is perhaps most natural to ask for geometric
realisations in Calabi-Yau varieties over $\Q$. For weight 3,
i.e.~dimension 2, this leads exactly to the singular K3 surfaces
that we started with.
In this case Question \ref{Question} was formulated independently by Mazur and van Straten.

At the point when this classification result for CM-forms was established, there were 35 forms known to correspond to singular K3 surfaces over $\Q$. These included all forms for fields of class number 1 and 2 through Kummer surfaces of products of elliptic curves. Further forms were realised in extremal elliptic surfaces (\cite{Sch}, \cite{BM}). In the meantime, N.~Elkies and the author found singular K3 surfaces over $\Q$ for all 65 known imaginary-quadratic fields with class group exponent 2 \cite{ES}. Since these singular K3 surfaces admit elliptic fibrations over $\Q$, all twists can be realised geometrically. 

\begin{Theorem}[Elkies, Sch\"utt]\label{Thm:ES}
Assuming ERH for odd real Dirichlet characters, every newform of weight 3 with rational coefficients is associated to a singular K3 surface over $\Q$.
\end{Theorem}

In case of CM newforms of weight 4 with rational coefficients, we ask for a corresponding (preferably rigid) Calabi-Yau threefold over $\Q$. A rigid Calabi-Yau partner of the newform of level $9$ from Table \ref{T:wt4} was established by Werner with assistance of van Geemen \cite{WvG}. For the forms of level 32, 49, and 256, Cynk and Meyer found non-rigid Calabi-Yau partners in \cite{CM}. Then Cynk and Hulek derived a generalised Kummer construction which realises all CM newforms of weight 4 and class number one \cite{CH}. Lately, this construction was combined with Weil restrictions by Cynk and the author \cite{CS}. This method produces non-rigid Calabi-Yau threefolds over $\Q$ for all CM newforms of weight $4$ except for the last newform in Table \ref{T:wt4} (which has class number $9$, i.e.~it is the only known imaginary-quadratic field of class group exponent $3$, but class number greater than $3$).

For higher dimensions
even less is known. Only recently, Cynk and Hulek exhibited a
general approach \cite{CH}. One particular construction involves
the $n$-fold product of the elliptic curve $E$ with CM by the ring of integers
in ${\Q(\sqrt{-3})}$. After dividing out
$E^n$ by a certain finite group action, Cynk and Hulek derive a
resolution of singularities $X$, which is Calabi-Yau and defined
over $\Q$. For $n$ odd, this has middle cohomology group
\[
\text{H}^n(X)=\text{H}^{n0}(X)\oplus\text{H}^{0n}(X)
\]

\pagestyle{empty}

which is two-dimensional. Its $L$-series comes from the Hecke character of $\Q(\sqrt{-3})$ in Lemma \ref{Lem:surj}. If $n$ is even, the corresponding statement can be established for the transcendental lattice $T(X)\subset\text{H}^n(X,\Z)$.

\pagestyle{headings}

\vspace{0.4cm}

\textbf{Acknowledgement:} I am indepted to K. Hulek for his
continuous interest and encouragement. Partial support by the DFG
Schwerpunkt 1094 ''Globale Methoden in der komplexen Geometrie" is
gratefully acknowledged. My thanks go also to the referee for helpful comments.

Part of the revising took place while I enjoyed the hospitality of
the Dipartimento di Matematica "Frederico Enriques" of Milano
University. Funding from the network Arithmetic Algebraic
Geometry, a Marie Curie Research Training Network, is gratefully
acknowledged. I particularly thank M. Bertolini and B. van Geemen.
The final version was prepared while I was funded by DFG
under grant Schu 2266/2-2.

\vspace{0.5cm}

Matthias Sch\"utt\\
Department of Mathematical Sciences\\
University of Copenhagen\\
Universitetspark 5\\
2100 Copenhagen\\
Denmark\\
{\tt mschuett@math.ku.dk}

\end{document}